\documentclass[11pt]{article}
\usepackage%[colorlinks]
{hyperref}

\author{M.~Kazarian}
\title{KP hierarchy for Hodge integrals}
\date{}

\usepackage{amssymb,amsmath}
\usepackage[matrix,arrow]{xy}
\usepackage{graphicx}

\def\C{{\mathbb C}}

\def\Z{{\mathbb Z}}

\def\ocM{{\overline{\cal M}}}

\def\cF{{\cal F}}
\def\const{{\rm const}}

\usepackage{theorem}

\newtheorem{theorem}{Theorem}[section]
\newtheorem{lemma}[theorem]{Lemma}

\newtheorem{proposition}[theorem]{Proposition}
\newtheorem{corollary}[theorem]{Corollary}

{\theorembodyfont{\rmfamily}
\newtheorem{remark}[theorem]{Remark}
\newtheorem{example}[theorem]{Example}

}

\let\l\lambda
\let\b\beta

\let\g\gamma
\let\f\varphi

\let\d\delta
\let\t\tau
\def\pd#1#2{\frac{\partial#1}{\partial#2}}
\let\D\partial

\def\ocM{{\overline{\mathcal M}}}
\let\L\Lambda

\def\ocM{{\overline{\cal M}}}
\def\gl{{\frak g\frak l}(\infty)}
\let\s\sigma
\let\^\wedge
\let\h\widehat
\def\ch{{\rm ch}}

\begin{document}
\maketitle

\begin{abstract}
Starting from the ELSV formula, we derive a number of new
equations on the generating functions for Hodge integrals over the
moduli space of complex curves. This gives a new simple and
uniform treatment of certain known results on Hodge integrals like
Witten's conjecture, Virasoro constrains, Faber's
$\l_g$-conjecture etc. Among other results we show that a properly
arranged generating function for Hodge integrals satisfies the
equations of the KP hierarchy.
\end{abstract}

\section{Introduction}

By Hodge integrals we mean intersection numbers of the form
$$\langle\l_j\t_{k_1}\dots\t_{k_n}\rangle=
 \int_{\ocM_{g,n}}\l_j\psi_1^{k_1}\dots\psi_n^{k_n},$$
where $\ocM_{g,n}$ is the moduli space of complex stable curves
with $n$ ordered marked points,~$\psi_i$ is the first Chern class
of the line bundle over $\ocM_{g,n}$ formed by the cotangent lines
at the $i$th marked point and~$\l_j$ is the $j$th Chern class of
the rank~$g$ Hodge vector bundle whose fibers are the spaces of
holomorphic one-forms. These numbers are well defined whenever the
equality $j+\sum k_i=3g-3+n(=\dim \ocM_{g,n})$ holds. They play an
important role in various problems related to the Gromov-Witten
theory.

There are several approaches to intersection theory on moduli
spaces. Among those approaches the one that seems to be the most
simple and the most straightforward is that based on the
Ekedahl-Lando-Shapiro-Vainshtein (ELSV) formula~\cite{ELSV}. This
formula expresses Hurwitz numbers enumerating ramified coverings
of the sphere as linear combinations of Hodge integrals. The ELSV
formula can be inverted in order to get some information about
Hodge integrals from known facts about Hurwitz numbers.

The ELSV formula was applied successfully in many papers, see,
e.g.~\cite{OP}, \cite{KL}, \cite{GJV2}, \cite{CLL}, \cite{SZ}. In
this note, we undertake a revision of the methods developed in
those papers. One of our main goals is to give a complete and
clear description of the relationship between various known
equations for Hurwitz numbers on one side and Hodge integrals on
the other side: equations of integrable hierarchies, the
cut-and-join equation, Virasoro constrains, $\l_g$-conjecture etc.

The method used in the present paper is close to that
of~\cite{KL}. The main difference is in a new change of variables
inverting the ELSV formula. The advantage of the new change is
that it induces an automorphism of the KP hierarchy. This permits
one to derive the whole hierarchy of PDE's at once for the
generating function of Hodge integrals, and, as a consequence, for
the Witten's potential participating in his conjecture.

The change of variables explored in the present paper is motivated
by that of Goulden-Jackson-Vakil used in~\cite{GJV2} in their
proof of the $\l_g$-conjecture and also in the paper~\cite{CLL}
devoted to the new derivation of the Virasoro constrains for
intersection numbers of $\psi$ classes. In both papers the GJV
change is done via the so called symmetrization operation. The
symmetrization is used as an intermediate step of computations and
is not used in the formulation of the final result. It is natural,
therefore, to try to skip the symmetrization operation and to
apply the GJV change directly to the original generating function.
The main difficulty appearing in this approach is technical: one
has to make a change of variables in a differential operator
containing infinitely many summands and involving infinitely many
variables. To overcome this difficulty, we apply here the
machinery of the boson-fermion correspondence. It allows one to
reduce the manipulation with differential operators in an infinite
dimensional space to those in just one variable and containing
finitely many terms. This makes all computations quite elementary
and free of combinatorial difficulties (no infinite sums are
involved).

I would like to thank my colleague S. Lando and all participants
of our joint seminar at the Independent University of Moscow where
the preliminary version of the presented theory was intensively
discussed.

\section{Main results}

Let us collect Hodge integrals into the following formal series in
an infinite set of formal commuting variables $u,T_0,T_1,\dots$:
\begin{equation}\label{preHodge}
\sum_{j,k_0,k_1,\dots}(-1)^j
 \langle \l_j\,\t_{0}^{k_0}\,\t_{1}^{k_1}\dots\rangle u^{2j}
 \frac{{T_0}^{k_0}}{k_0!}\frac{{T_1}^{k_1}}{k_1!}\dots,
\end{equation}
where the summation is taken over all possible monomials in the
symbols~$\t_i$ and over all possible values~$j\ge0$. Denote by
$G(u;q_1,q_2,\dots)$ the series obtained from~\eqref{preHodge} by
the following linear substitution of variables
\begin{align*}
 T_0&=q_1,\\
 T_1&=u^2\,q_1 + 2\,u\,q_2 + q_3,\\
 T_2&=u^4\,q_1 + 6\,u^3\,q_2 + 12\,u^2\,q_3 + 10\,u\,q_4 + 3\,q_5,\\
 T_3&=u^6\,q_1 + 14\,u^5\,q_2 + 61\,u^4\,q_3 + 124\,u^3\,q_4
   + 131\,u^2\,q_5 + 70\,u\,q_6 + 15\,q_7,\\
   &\qquad\dots.
\end{align*}

In general, the linear combination $T_k$ of the variables $q_i$ is
defined as follows (the meaning of this change will be explained
in the next section). Consider the sequence of polynomials
$\f_k(u,z)$, $k=0,1,2,\dots$, defined by
$$\f_0(u,z)=z,\quad \f_{k+1}(u,z)=D\,\f_{k}(u,z)=D^{k+1}\,\f_0(u,z),\quad
 \hbox{where }D=(u+z)^2z\,\pd{}{z}:$$
$$\f_0=z,\quad \f_1=u^2z+2uz^2+z^3,\quad
 \f_2=u^4z+6u^3z^2+12u^2z^3+10uz^4+3z^5,\quad\dots.$$
Then $T_k$ is obtained form $\f_k$ by replacing $z^m$ by $q_m$ in
each monomial. Equivalently, $T_k$ is given by the following
recursive equation
\begin{equation}\label{D}
T_{k+1}=\sum_{m\ge1}m\,(u^2\,q_{m}+2\,u\,q_{m+1}+q_{m+2})\,\pd{}{q_m}T_k,
\end{equation}
By construction, $T_k$ is a linear combination of variables~$q_s$,
besides, the variable with the maximal index~$s=2k+1$ enters with
the coefficient~$(2k-1)!!$, and the coefficients of the variables
with smaller indices contain positive powers of the parameter~$u$.

Remark that the result of substitution $u=0$ to $G$ only depends
on variables $q_k$ with odd~$k$ and it turns into the Witten's
potential~$F$ for the intersection numbers of~$\psi$ classes after
rescaling $q_{2d+1}=\frac{t_d}{(2d-1)!!}$.

\begin{theorem}\label{th2}
The series~$G$ is a solution of the KP hierarchy with respect to
the variables~$q_i$ {\rm(}identically in~$u${\rm)}.
\end{theorem}

Witten's conjecture (now Kontsevich's theorem, see~\cite{W},
\cite{Ko}) claims that~$F$ is a solution of the KdV hierarchy.
This statement is an obvious specialization of the previous
theorem. Indeed, the equations of the KdV hierarchy are obtained
from the equations of the KP hierarchy by an additional
requirement that the function is independent of even variables.

\medskip
As it was shown by C. Faber~\cite{F} (based on earlier result of
Mumford~\cite{Mu}, see Section~\ref{secHodge} below), the
computation of Hodge integrals can be reduced to the computation
of the intersection indices of $\psi$ classes. In other words, all
coefficients of the series~$G$ are determined by the coefficients
of~$F$. Therefore, one can try to derive Theorem~\ref{th2} from
the statement of Witten's conjecture. However, our direct
arguments are based on the application of the ELSV formula
relating Hodge integrals to the Hurwitz numbers. The Hurwitz
numbers participating in this formula are discussed in the next
section. Here we only remark that they are relatively simple
combinatorial objects, in particular, the generating series
$H(\b;p_1,p_2,\dots)$ for these numbers can be given by the
following explicit closed formula
\begin{equation}\label{eHM0}
 e^H=e^{\b M_0}e^{p_1},
\end{equation}
where $M_0$ is the so-called \emph{cut-and-join operator},
$$M_0=\frac12\sum_{i,j}\Bigl((i+j)p_ip_j\pd{}{p_{i+j}}+
 i\,j\,p_{i+j}\pd{^2}{p_i\partial p_j}\Bigr).$$
%~(\cite{Ok})
%\begin{equation}\label{eH}
%e^H=\sum_\l s_\l(1,0,\dots)e^{f_2(\l)\b}s_\l(p_1,p_2,\dots),
%\end{equation}
%where the summation is taken over the set of all Young diagrams
%(partitions) $\l=(\l_1,\l_2,\dots)$, $\l_1\ge\l_2\ge\dots\ge0$;
%$s_\l$ is the Schur function corresponding to the partition~$\l$,
%and
%\begin{equation}\label{f2}
%f_2(\l)=\frac12\sum_{i=1}^\infty\Bigl((\l_1-i+\tfrac12)^2
% -(-i+\tfrac12)^2\Bigr).
%\end{equation}

As we shall see in Section~\ref{KPsec}, the very existence of a
formula like~\eqref{eHM0} implies immediately the following
statement.

\begin{theorem}{\rm(cf.~\cite{Ok,KL})}\label{th3}
The generating function~$H$ for Hurwitz numbers is a solution of
the KP hierarchy with respect to the variables $p_i$
{\rm(}identically in the formal parameter~$\b${\rm)}.
\end{theorem}

The ELSV formula expresses the Hurwitz numbers (the coefficients
of~$H$) in terms of the Hodge integrals (the coefficients of~$G$).
By formal manipulations this formula can be reduced to the
following one.

Consider two variables $x$ and $z$ related to one another by the
changes of variables
\begin{equation}\label{xz}\begin{aligned}
 x&=\frac{z}{1+\b\,z}e^{-\frac{\b\,z}{1+\b\,z}}=z - 2\,\b\,z^2
 + \frac72\,\b^2\,z^3
 - \frac{17}3\,\b^3\,z^4 +\dots,\\
 z&=\sum_{b\ge1}\frac{b^b}{b!}\b^{b-1}x^b=x+2\,\b\,x^2
  +\frac92\,\b^2\,x^3+\frac{32}3\,\b^3\,x^4+\dots.
 \end{aligned}\end{equation}
The fact that these changes of variables are inverse to one
another follows from the Lagrange inversion theorem,
see~\cite{GJ,GJV1}. These changes provide a linear isomorphism
(depending on the parameter~$\b$) of the spaces of formal power
series in the variables~$x$ and~$z$. If we identify the linear
span of the variables~$p_i$ with the space of formal power series
in~$x$ and the linear span of the variables~$q_i$ with the space
of series in~$z$ by means of the correspondence
\begin{equation}\label{pxqz}
p_b\leftrightarrow x^b,\qquad
 q_k\leftrightarrow z^k,
\end{equation}
then the isomorphism above provides a linear change of variables
(depending on the parameter~$\b$) between variables~$p_i$
and~$q_j$. More explicitly, this change is given by
\begin{equation}\label{ch}
p_b= \sum_{k\ge b}c^b_k\,\b^{k-b}\,q_k,
\end{equation}
where the rational coefficients~$c^b_k$ are determined by the
expansion
$$x^b=\sum_{k\ge b}c^b_k\b^{k-b}\,z^k.$$

Let us set also
\begin{equation}\label{H01H02}
H_{0,1}=\sum_{b=1}^\infty \frac{b^{b-2}}{b!} p_b \b^{b-1},
 \qquad H_{0,2}=\frac12\sum_{b_1,b_2=1}^\infty
 \frac{b_1^{b_1}b_2^{b_2}}{(b_1+b_2)b_1!b_2!}
 p_{b_1}p_{b_2}\b^{b_1+b_2}.
\end{equation}

\begin{theorem}\label{th4}
The change~\eqref{ch} in $H-H_{0,1}-H_{0,2}$ leads to the
series~$G$ of Theorem~{\rm\ref{th2}}, up to the
rescaling~$q_k\mapsto\b^{\frac43k}q_k$ and $u\mapsto\b^{1/3}$,
where~$H=H(\b;p_1,p_2,\dots)$ is the generating function for
Hurwitz numbers:
$$(H-H_{0,1}-H_{0,2})|_{p=p(\b;q)}=G(\b^{\frac13};
 \b^{\frac43}q_1,\b^{\frac83}q_2,\b^{\frac{12}3}q_3,\dots).$$
\end{theorem}

\begin{remark}
The action of the change~\eqref{ch} described by Theorem~\ref{th4}
can be characterized as follows. Consider the plane of
coordinates~$(m,B)$ and mark all points of this plane
corresponding to non-trivial terms of the form
$\const\,p_{b_1}\dots p_{b_n}\b^m$ in the series~$H$, where
$B=\sum b_i$ (see the picture).
\begin{figure}[ht]
\centerline{\includegraphics[scale=0.75]{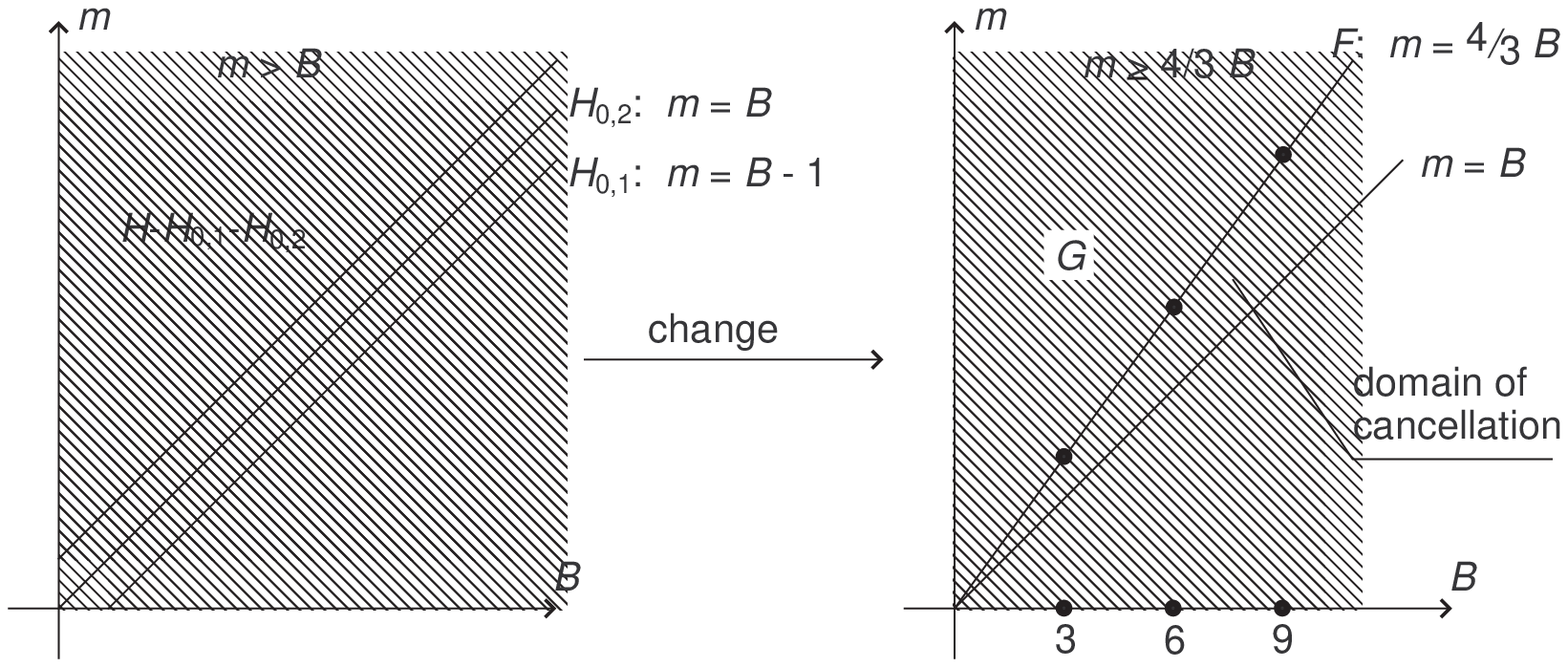}}
\end{figure}%
There are no marked points below the line $m=B-1$ . The points
lying on the lines $m=B-1$ and $m=B$ correspond to the
contributions of $H_{0,1}$ and $H_{0,2}$ to~$H$, respectively. The
points, corresponding to the contribution of the remaining terms
are situated above the diagonal $m=B$. The change~\eqref{ch}
determines a linear transformation of the coefficients in~$H$.
Furthermore, a coefficient of the original series may contribute
to some coefficient of the resulting series only if the
corresponding points lie on the same diagonal $m-B=\const$ and the
point corresponding to the original series has smaller
coordinates~$m$ and~$B$. After the application of this change to
the series $H-H_{0,1}-H_{0,2}$, we get, after some `magic
cancellations', a series having no nontrivial coefficients in the
domain $m<\frac43 B$. The terms of the resulting series lying on
the line $m=\frac43 B$ form the generating series~$F$ for
intersection numbers of~$\psi$ classes and the other terms
correspond to other summands of the series~$G$ of
Theorem~\ref{th2} (after the rescaling $q_i\mapsto
\b^{-\frac43i}q_i$, $\b^{\frac13}\mapsto u$).
\end{remark}

Consider now an arbitrary invertible formal series $x(z)$ and
associate to this series a linear change of variables $p\mapsto
p(q)$ given by
\begin{equation}
 p_b=\sum_{k\ge b} c_k^b q_k,
 \label{gench}\end{equation}
whose coefficients
$c_k^b$ are determined by the expansion $x^b=\sum_{k\ge b}c_k^b
z^k$.

\begin{theorem}\label{th5}
There is a quadratic function $Q(p_1,p_2,\dots)$ such that the
transformation sending an arbitrary series $\Phi(p_1,p_2,\dots)$
to the series $\Psi(q_1,q_2,\dots)=(\Phi+Q)|_{p\to p(q)}$ is an
authomorphism of the KP hierarchy: it sends solutions to
solutions.

In the case when
$x(z)=\frac{z}{1+\b\,z}e^{-\frac{\b\,z}{1+\b\,z}}$ this quadratic
function is $Q=-H_{0,2}$.
\end{theorem}

Theorem~\ref{th2} is an obvious corollary of
Theorems~\ref{th3},~\ref{th4}, and~\ref{th5}. Indeed, the
transformation of Theorem~\ref{th4} differs from that of
Theorem~\ref{th5} by linear terms $H_{0,1}$ that do not affect the
equations of the KP hierarchy. The rescaling
$q_k\mapsto\b^{\frac43k}q_k$ is also an automorphism of the
hierarchy since its equations are quasihomogeneous. Therefore, the
validity of these equations for the series~$G$ is equivalent to
their validity for~$H$, which is guaranteed, in turn, by the
assertion of Theorems~\ref{th3}.

\begin{remark}
The definition for the change~\eqref{ch} looks unmotivated.
Implicitly, the same change was used in~\cite{GJV1},~\cite{GJV2}.
The only motivation that we can provide here is that `it works'.
In fact, there is a freedom in the choice of a change inverting
the ELSV formula. One of the possible changes was used
in~\cite{KL}. That change did not preserve the whole KP hierarchy
but was sufficient to derive the KdV equation for the Witten's
potential~$F$.
\end{remark}

\section{The ELSV formula}

This section is devoted to the proof of Theorem~\ref{th4}.

Consider a ramified covering of the sphere~$S^2=\C P^1$ by a
smooth surface of genus~$g$ such that the point~$\infty\in \C P^1$
has~$n$ marked preimages of multiplicities~$b_1,\dots,b_n$ and all
other critical points are simple (with ramification of the second
order each) and have pairwise different critical values. The
number~$m$ of simple critical values is determined by the
Riemann-Hurwitz formula:
$$m=2g-2+n+\sum_{i=1}^n b_i.$$
The Hurwitz number $h_{g;b_1,\dots,b_n}$ is defined as the number
of such coverings with a fixed position of the critical values;
the coverings are counted with their weights inverse to the order
of the automorphism group of the covering. The celebrated ELSV
formula~\cite{ELSV} expresses these numbers via Hodge integrals:
$$\frac{h_{g;b_1,\dots,b_n}}{m!}=
 \prod_{i=1}^n\frac{b_i^{b_i}}{b_i!}\int_{\ocM_{g,n}}
 \frac{1-\l_1+\l_2-\dots\pm\l_g}{\prod_{i=1}^n(1-b_i\psi_i)}.$$
It is convenient to pack Hurwitz numbers into the generating
series
$$H(\b;p_1,p_2,\dots)=\sum_{n\ge1}\frac1{n!}\sum_{g,b_1,\dots,b_n}
 h_{g;b_1,\dots,b_n}\frac{\b^m}{m!}p_{b_1}\dots p_{b_n}.$$
The ELSV formula is not applicable to coverings of genus~$0$
with~$1$ or~$2$ preimages at infinity (since the corresponding
moduli spaces $\ocM_{0,1}$ and $\ocM_{0,2}$ do not exist). The
summands $H_{0,1}$ and $H_{0,2}$, respectively, corresponding to
such coverings are presented in Eq.~\eqref{H01H02} of the previous
section. Let us represent the remaining summands in the form
$H-H_{0,1}-H_{0,2}=\sum_{n\ge1}\frac{1}{n!}H_n$, where $H_n$
contains all terms corresponding to the coverings with exactly~$n$
preimages at infinity. Then, since
%$\dim\ocM_{g,n}=3g-3+n$ ¨
$$m=2g-2+n+\sum b_i
 =\sum(b_i+\tfrac13)+\tfrac23(3g-3+n)
 =\sum(b_i+\tfrac13)+\tfrac23\dim\ocM_{g,n},$$
we get
\begin{align*}
 H_n&=\sum_{g,b_i,\dots,b_n}\prod_{i=1}^n
  \frac{b_i^{b_1}\,\b^{b_i+\frac13}}{b_i!}
  \int_{\ocM_{g,n}}\frac{1-\b^{\frac23}\l_1+\b^{\frac43}\l_2-\dots}
  {\prod_{i=1}^n(1-b_i\b^{\frac23}\psi_i)}p_{b_1}\dots p_{b_n}\\
  &=\left\langle
  (1-\b^{\frac23}\l_1+\b^{\frac43}\l_2-\dots)
  \prod_{i=1}^n\sum_{b\ge1}
  \frac{b^{b}}{b!}\frac{\b^{b+\frac13}p_b}
  {(1-b\b^{\frac23}\psi_i)}
  \right\rangle\\
  &=\left\langle
  (1-\b^{\frac23}\l_1+\b^{\frac43}\l_2-\dots)
  \prod_{i=1}^n\sum_{d\ge0}T_d\psi_i^d
  \right\rangle,
\end{align*}
where
$$T_d=\sum_{b\ge1}\frac{b^{b+d}}{b!}\b^{b+\frac13+\frac23d}p_b,$$
and where for each monomial in the classes $\psi_i$ and $\l_j$ we
denote by angular brackets its integral over the space
$\ocM_{g,n}$, $g\ge0$, whose dimension is equal to the degree of
the monomial. Setting $u=\b^{\frac13}$ we can rewrite the ELSV
formula in the following form (a similar form of the ELSV formula
was observed in~\cite{GJV1})
$$ H-H_{0,1}-H_{0,2}=\sum_{j,k_0,k_1,\dots}(-1)^j
 \langle \l_j\,\t_{0}^{k_0}\,\t_{1}^{k_1}\dots\rangle
 u^{2j}\frac{{T_0}^{k_0}}{k_0!}\frac{{T_1}^{k_1}}{k_1!}\dots.
$$
It follows that Theorem~\ref{th4} is a corollary of the following
lemma.

\begin{lemma}\label{3.1}
The defined above series $T_d$ in the variables~$p_i$ turn under
the change~\eqref{ch} into polynomials in variables $q_i$
participating in the statement of Theorem~{\rm\ref{th2}}, up to
the rescaling $q_i\mapsto \b^{\frac43i}q_i$,
$u\mapsto\b^{\frac13}$.
\end{lemma}

\emph{Proof}. Under the correspondence~\eqref{pxqz}, the
series~$T_d$ corresponds to the following function in~$x$:
$$T_d\leftrightarrow
 \b^{\frac43+\frac23d}\sum_{b\ge1}\frac{b^{b+d}}{b!}\b^{b-1}x^b
 =(\b^{\frac23}D)^d\b^{\frac43}z(x),$$
where $D=x\pd{}{x}$ and where the series
$z(x)=\sum_{b\ge1}\frac{b^{b}}{b!}\b^{b-1}x^b$ is inverse to the
series $x(z)=\frac{z}{1+\b z}e^{-\frac{\b z}{1+\b z}}$. It is easy
to check that the differential operator~$D=x\pd{}{x}$ takes in
terms of the coordinate~$z$ the form
$$D=(1+\b\,z)^2z\pd{}{z}.$$
Under the correspondence~\eqref{pxqz} the function
$\b^{\frac43}z\leftrightarrow T_0$ turns into~$\b^{\frac43}q_1$,
and the operator~$\b^{\frac23}D$ turns, respectively,~into
\begin{align*}
\b^{\frac23}D&=\sum_{m\ge1}m\,(\b^{\frac23}\,q_{m}+2\,\b^{\frac53}\,q_{m+1}
 +\b^{\frac83}\,q_{m+2})\,\pd{}{q_m}\\
 &=\sum_{m\ge1}m\,(\b^{\frac23+\frac43m}q_{m}
 +2\,\b^{\frac13+\frac43(m+1)}\,q_{m+1}
 +\b^{\frac43(m+2)}\,q_{m+2})
 \frac1{\b^{\frac43m}}\,\pd{}{q_m}.
\end{align*}
After rescaling $\b^{\frac43m}q_m\mapsto q_m$,
$\b^{\frac13}\mapsto u$ this operator coincides with the operator
participating in the recursive relation~\eqref{D}. This proves
Lemma~\ref{3.1}, and hence, Theorem~\ref{th4}.

\section{KP hierarchy}\label{KPsec}

In this section we recall several well known facts about the KP
hierarchy that are sufficient for the proof of Theorems~\ref{th3}
and~\ref{th5}. For a more detailed exposition of the theory we
refer to the papers~\cite{DKJM}, \cite{MJD}, and to
Section~\ref{sato} of the present paper. The KP
(Kadomtsev-Petviashvili) hierarchy is a particular system of
partial differential equations on the unknown function (power
series)~$F$ in infinite set of variables $p_1,p_2,\dots$. Here are
several first equations of the hierarchy
 {%\small
 \begin{equation}
\begin{aligned}\label{KP}
 F_{2,2}&=-\frac12\,F_{1,1}^2+F_{3,1}-\frac1{12}\,F_{1,1,1,1},\\
 F_{3,2}&=-F_{1,1}F_{2,1}+F_{4,1}-\frac16F_{2,1,1,1},\\
 F_{4,2}&=-\frac12\,F_{2,1}^2-F_{1,1}F_{3,1}+F_{5,1}+\frac18\,F_{1,1,1}^2\\
  &\qquad\qquad
     +\frac1{12}\,F_{1,1}F_{1,1,1,1}-\frac14\,F_{3,1,1,1}+\frac1{120}\,F_{1,1,1,1,1,1},\\
 F_{3,3}&=\frac13\,F_{1,1}^3-F_{2,1}^2-F_{1,1}F_{3,1}+F_{5,1}
     +\frac14\,F_{1,1,1}^2\\
  &\qquad\qquad
     +\frac13\,F_{1,1}F_{1,1,1,1}-\frac13\,F_{3,1,1,1}+\frac1{45}\,F_{1,1,1,1,1,1}.
\end{aligned}
\end{equation}}

The exponent $\t=e^F$ of any solution is called the
\emph{$\t$-function} of the hierarchy. It is known that the space
of solutions (or the space of $\t$-functions) is homogeneous:
there is a Lie algebra $\h\gl$ acting on the space of solutions
and the action of the corresponding Lie group is transitive. In
other words, any solution can be obtained from any other solution
(say, from the solution $\t=1$) by the action of an appropriate
transformation from this group. The actual definition of the Lie
algebra $\h\gl$ and of its action is given in Section~\ref{sato}.
Here we only present several sample operators (acting on
$\t$-functions) belonging to this algebra.

\begin{example}
The scalar operator (the operator of multiplication by a constant)
belongs to $\h\gl$. It follows that the $\t$-function is defined
up to a multiplicative constant. This constant is usually chosen
in such a way that $\t(0)=1$. In this case the logarithm
$F=\log(\t)$ is a correctly defined power series if the
function~$\t$ is.
\end{example}

\begin{example}
Set
$$a_k=\begin{cases}
 p_k&k>0\\
 0&k=0\\
 (-k)\pd{}{p_{-k}}&k<0
 \end{cases}$$
Then for all~$k$ the operator $a_k$ belongs to $\h\gl$. For
positive~$k$ this means that addition of a linear function
preserves the space of solutions of the hierarchy. For
negative~$k$ this means that a shift of arguments also preserves
the hierarchy. Both assertions are obvious. Indeed, the equations
of the hierarchy have constant coefficients and the partial
derivatives have order at least two.
\end{example}

\begin{example}
The following operator also belongs to $\h\gl$ for any integer
$m\ne0$:
$$\L_m=\frac12\sum_{i=-\infty}^{\infty}a_ia_{m-i}.$$
The right hand side is well defined for any~$m\ne0$ since~$a_i$
and~$a_{m-i}$ commute. For $m=0$ the formula should be corrected:
$$\L_0=\sum_{i=1}^\infty
 a_{i}a_{-i}=\sum_{i=1}^\infty i\,p_i\pd{}{p_i}.$$
For positive~$m$ these operators have the following form
\begin{equation}
\L_m=\sum_{i\ge 1}
 i\,p_{i+m}\pd{}{p_i}+\frac12\sum_{j=1}^{m-2}p_jp_{m-1-j},
 \qquad m>0.
 \label{Lm}\end{equation}
For negative $m$ the operators $\L_m$ involve second order partial
derivatives. The operators $\L_{-2m}$, $m\ge-1$, participate in
the Virasoro equations for the Witten's potential of intersection
indices of~$\psi$ classes (see Section~\ref{secVir}).
%In what follows we
%shall use the following operators
%\begin{align*}
%\L_{1}&=\sum_{i\ge1}i\,p_{i+1}\pd{}{p_i},\\
%\L_{2}&=\frac{p_1^2}2+\sum_{i\ge1}i\,p_{i+2}\pd{}{p_i}.
%\end{align*}
\end{example}

\begin{example}\label{eH1}
The following operators belong to~$\h\gl$ for any integer~$m$,
$$M_m=\frac16\sum_{i,j=-\infty}^\infty {:}a_ia_ja_{m-i-j}{:},$$
where we use notation ${:}a_{i_1}\dots a_{i_k}{:}=a_{\s(i_1)}\dots
a_{\s(i_k)}$, where $\s$ is a permutation of the indices
$i_1,\dots,i_k$ such that $\s(i_1)\ge\dots\ge \s(i_k)$. The
operator $M_0$ is called also the \emph{cut-and-join} operator
$$
 M_0=\frac12\sum_{i,j=1}^\infty\Bigl((i+j)\,p_ip_j\pd{}{p_{i+j}}+
 i\,j\,p_{i+j}\pd{^2}{p_i\D p_j}\Bigr).$$
It is known (see, e.g. \cite{GJ1}) that the exponent of the
generating function for Hurwitz numbers satisfies the cut-and-join
equation
\begin{equation}\label{caj}
\pd{e^H}{\b}=M_0e^H.
\end{equation}
Since for $\b=0$ the function $H|_{\b=0}=p_1$ satisfies the KP
hierarchy by trivial reason and since the infinitesimal
transformation of~$H$ with a change of the parameter~$\b$ belongs
to~$\h\gl$, we conclude that~$H$ satisfies the KP hierarchy for
all parameter values~$\b$. This proves, in particular,
Theorem~\ref{th3}.

More explicitly, the solution of the cut-and-join equation is
given by
$$e^H=e^{\b\,M_0}e^{p_1}.$$
This expression shows that the $\t$-function $e^H$ is obtained
from the trivial one $1=e^0$ by the action of the composition of
the operators $e^{a_1}$ and $e^{\b M_0}$ both of which belong to
the Lie group of the algebra $\h\gl$.

% In order to compute the
%exponent of~$M_0$ explicitly, remark that the Schur functions form
%an eigenbasis for~$M_0$ and the corresponding eigenvalues are
%given by~\eqref{f2}. Therefore, the formula above for~$e^H$ is
%equivalent to~\eqref{eH}, see~\cite{KL}.
\end{example}

We are ready now to fulfil the computations leading to the proof
of Theorem~\ref{th5}. An infinitesimal version of the change
$z\mapsto x(z)$ is a vector field on the line of the
coordinate~$z$ which can be written as a linear combination of the
basic fields $z^{m+1}\pd{}{z}$, $m\ge0$. These fields can be
viewed as linear transformations of the space of power series in
$z$:
$$z^{m+1}\pd{}{z}:z^i\mapsto i\,z^{i+m}.$$
Under the correspondence~\eqref{pxqz} this operator sending $p_i$
to $i\,p_{i+m}$ can be written as a linear differential operator
$\sum i\,p_{i+m}\pd{}{p_i}$ which is nothing but the
`differential' part of the operator~$\L_m$~\eqref{Lm}. Integrating
such infinitesimal transformations we obtain a global linear
change of variables of the form~\eqref{gench}. The `polynomial'
part of~$\L_m$ is a quadratic form. Since the space of quadratic
functions is invariant under linear changes of variables, the
first statement of Theorem~\ref{th5} follows.

To prove the second statement, we observe that the change
$x(z)=\frac{z}{1+\b z}e^{-\frac{\b z}{1+\b z}}$ can be obtained
from the identity as the time $\b$ flow of the non-autonomous
vector field $-(2\,z+\b\,z^2)\,z\pd{}{z}$. This assertion is a
reformulation of the following easily verified identity
$$\pd{x(z)}{\b}=-(2\,z+\b\,z^2)\,z\pd{x(z)}{z}.$$
The field~$-(2\,z+\b\,z^2)\,z\pd{}{z}$ corresponds to the operator
$-(2\,\L_{1}+\b\,\L_{2})\in\h\gl$. It remains to check that the
quadratic function $Q=-H_{0,2}$ agrees with the `polynomial' part
of this operator.

Consider an arbitrary series $\Phi(p_1,p_2,\dots)$ and set
$Z(\b;q)=\exp(\Phi-H_{0,2})|_{p\to p(\b;q)}$, where the linear
change of coordinates $p\to p(\b;q)$ corresponds to our choice
of~$x(z)$.

\begin{lemma}\label{inf}
The series $Z$ is subject to the differential equation
$$\pd{Z}{\b}=-(2\,\L_{1}+\b\,\L_{2})Z.$$
\end{lemma}

\emph{Proof}. Differentiating the function~$Z=e^{\Phi-H_{0,2}}$ we
get
\begin{equation}\label{dZ}
\pd{Z}{\b}=\Bigl(
 -\pd{H_{0,2}}{\b}\Bigm|_{p=\const}
 +\sum_{b\ge1}\pd{p_b}{\b}\pd{}{p_b}\Bigr)Z.
\end{equation}

We wish to rewrite the right hand side in $q$-coordinates. The
arguments above show that the second summand is the `differential'
part of the operator $-(2\,\L_{1}+\b\,\L_{2})$:
$$\sum_{b\ge1}\pd{p_b}{\b}\pd{}{p_b}=
 -\sum_{i\ge1}i\,(2\,q_{i+1}+\b\,q_{i+2})\pd{}{q_i}.$$
(The reader can check this equality by straightforward
 computations.) For the first summand, we have
$$-\pd{H_{0,2}}{\b}\Bigm|_{p=\const}=
 -\frac12\sum_{i,j\ge1}\frac{i^i}{i!}\frac{j^j}{j!}
 \b^{i+j-1}p_ip_j
 =-\frac\b2\,\Bigl(\sum_{i\ge1}\frac{i^i}{i!}\b^{i-1}p_i\Bigr)^2
 =-\frac\b2q_1^2,$$
which coincides with the `polynomial' part of the
operator~$-(\b\,\L_{2}+2\,\L_{1})$. The lemma is proved.

The lemma shows that $Z$ is a $\tau$-function for any parameter
value~$\b$ if it is for the initial parameter value $\b=0$. This
completes the proof of Theorem~\ref{th5}.

\section{Cut-and-join and the Virasoro constrains}\label{secVir}

Consider the cut-and-join equation~\eqref{caj} for the generating
function of Hurwitz numbers. The change of Theorem~\ref{th4}
transforming the series~$H$ to the series $G$ of Theorem~\ref{th2}
acts on differential equations as well. The cut-and-join equation
is transformed under this change into the following one.

\begin{theorem}\label{newcaj}
The exponent of the function~$G$ of Theorem~{\rm\ref{th2}} is
subject to the differential equation
\begin{multline*}
\frac13u^{-2}\pd{e^G}{u}=
 \Bigl(M_0+4u^{-1}M_1+6u^{-2}M_2+4u^{-3}M_3+u^{-4}M_4\\
 -\frac43u^{-3}\L_0-u^{-4}\L_1
 +\frac14u^{-2}a_2+\frac13u^{-3}a_3+\frac18u^{-4}a_4\Bigr)e^G,
\end{multline*}
where the operators $a_m$, $\L_m$, and $M_m$ are defined in the
previous section.
\end{theorem}

One of the possible proofs is the direct substitution. The reader
can try to compute himself the coefficients of the equation
starting from the definition of the transformation sending~$H$
to~$G$. More elementary computations leading to the same equation
are explained in Section~\ref{transf}.

\medskip
Each term in the equation of the theorem has degree greater than
or equal to $-4$ with respect to the variable~$u$. Since the
series~$F=G|_{u=0}$ only depends on odd variables, we obtain from
the explicit form of~$M_4$ that the variables $p_{2m}$ enter the
coefficient of $u^{-4}$ at most linearly. Extracting the
coefficient of $p_{2m+4}u^{-4}$ in the equation we get the
equalities
$$(2m+3)\pd{e^F}{p_{2m+3}}=
   \bigl(\L_{-2m}+\frac18\d_{m,0}\Bigr)e^F,\qquad m\ge-1.$$
These equations known as \emph{Virasoro constrains} for the
generating function~$F$ of the intersection numbers
of~$\psi$-classes form an equivalent reformulation of Witten's
conjecture. The derivation of the Virasoro constrains presented
above is parallel to that from~\cite{CLL}. I hope, however, that
in the presented form the computations look more clear and
natural.

\section{Boson-fermion correspondence}\label{sato}

In this section we review the basics of the boson-fermion
correspondence with application to Sato Grassmannian and the KP
hierarchy. The basic references are~\cite{DKJM} and~\cite{MJD}.

Consider the space $\C[[p_1,p_2,\dots]]$ of formal power series
and consider the additive basis in this space formed by Schur
functions
$$\C[[p_1,p_2,\dots]]\simeq\overline{\bigoplus_\l \C\,s_\l(p)}.$$
Schur functions are certain polynomials labelled by partitions
(Young diagrams). One of the possible their definitions is given
below. Here are several of them
\begin{gather*}
 s_0=1,\qquad s_1=p_1,\qquad s_2=\frac12(p_1^2+p_2),\qquad
 s_3=\frac16(p_1^3+3p_1p_2+2p_3),\\
 s_{1,1}=\frac12(p_1^2-p_2),
 \qquad s_{2,1}=\frac13(p_1^3-p_3), \qquad
 s_{1,1,1}=\frac16(p_1^3-3p_1p_2+2p_3).
\end{gather*}

There is another space with the basis naturally labelled by Young
diagrams. Namely, consider first an auxiliary space
$V=\C[z^{-1}][[z]]$ of formal Laurent series. The
\emph{semi-infinite wedge space}~$\L^{\frac\infty2}V$ is, by
definition, (the completion of) the vector space whose basic
vectors are semi-infinite formal wedge products of the form
$$v_\l=z^{k_1}\^z^{k_2}\^\dots,\qquad k_i=\l_i-i.$$
The sequences $(k_1,k_2,\dots)$ appearing in these products are
just arbitrary strictly decreasing sequences of integers such that
$k_i=-i$ for sufficiently large~$i$. The elements of the
space~$\L^{\frac\infty2}V$ can be represented as linear
combinations of infinite wedge products of the
form~$\f_1(z)\^\f_2(z)\^\dots$, $\f_i\in V$, such that
$\f_i=z^{-i}+{}$(terms of higher order in~$z$) for sufficiently
large~$i$. Using polylinearity and skew-symmetry of the wedge
product one can represent such a wedge product as a (possibly
infinite) linear combination of basic ones.

The \emph{boson-fermion correspondence} is the coordinate-wise
isomorphism of vector spaces
$$\C[p_1,p_2,\dots]\simeq\L^{\frac\infty2}V,
 \qquad s_\l\leftrightarrow v_\l.$$
The spaces on the left and the right hand sides of the isomorphism
are called \emph{bosonic} and \emph{fermionic Fock spaces} (of
zero charge), respectively. The vector
$v_\varnothing=z^{-1}\^z^{-2}\^z^{-3}\^\dots$ corresponding to the
function $s_\varnothing=1$ is called the \emph{vacuum vector}.

The geometric viewpoint to the theory of KP hierarchy is
formulated as follows.

\begin{theorem}
The function~$\t\in\C[[p_1,p_2,\dots]]$ is the exponent of a
solution of the KP hierarchy if and only if its image under
boson-fermion correspondence can be represented by a decomposable
wedge product
$$\t\leftrightarrow \f_1(z)\^\f_2(z)\^\dots.$$
\end{theorem}

Decomposable wedge products are characterized uniquely up to a
multiplicative constant by the linear span of the vectors~$\f_i$.
Therefore, the theorem has geometric reformulation that
$\t$-functions form the cone over the Grassmannian
$G_{\frac\infty2}(V)$ of half-infinite subspaces Pl\"ucker
embedded to the projective space $P\L^{\frac\infty2}(V)$. It is
known in algebraic geometry that the Pl\"ucker embedding of the
Grassmannian is given by quadratic equations. These algebraic
equations on the Taylor coefficients of the series $\t$ are called
\emph{Hirota bilinear equations}. They can be represented in a
form of partial differential equations on~$\t$. These equations
rewritten in terms of the logarithm $F=\log(\t)$ are exactly the
equations of the KP hierarchy.

\begin{example}
The wedge product $(z^{-1}+z^{2})\^z^{-2}\^z^{-3}\^\dots=v_0+v_3$
corresponds to the function
$$s_0+s_3=1+\frac16(p_1^3+3p_1p_2+2p_3).$$
It follows that its logarithm
$\log(1+\frac16(p_1^3+3p_1p_2+2p_3))$ satisfies the equations of
the KP hierarchy. The reader may check this fact by substituting
to the first of the equations~\eqref{KP}. These computations being
elementary occupy several pages and are quite laborious when made
by hand.
\end{example}

The Grassmannian is a homogeneous space; every subspace can be
obtained from any other by a linear transformation. The infinite
dimension of the spaces under consideration implies some
additional phenomena that we describe now. Denote by $\gl$ the Lie
algebra of differential operators in one variable~$z$ with Laurent
coefficients (completed in a suitable way that we do not discuss
here). Every element of this algebra can be treated as a linear
operator acting on the space~$V$. This operator is represented in
the basis~$z^i$, $i\in\Z$, by an infinite matrix $a_{i,j}$,
$z^j\mapsto \sum_i a_{ij} z^i$. This algebra is graded by the
agreement $\deg(z^i(\D/\D z)^j)=i-j$. For the operators of
degree~$k$ all non-zero components~$a_{i,j}$ are situated on the
diagonal $i-j=k$, moreover, the component~$a_{i,i-k}$ has
polynomial dependance on~$i$.

To every operator $A\in\gl$ we associate the operator $\h A$
acting on the fermionic space by the following rule. If the matrix
of~$A$ has no non-trivial diagonal elements, then the action of
$\h A$ is determined by the Leibnitz rule:
\begin{align*}
\widehat A(z^{k_1}\^z^{k_2}\^\dots)&
 =A(z^{k_1})\^z^{k_2}\^z^{k_3}\^\dots+\\
 &\qquad z^{k_1}\^A(z^{k_2})\^z^{k_3}\^\dots\\
 &\qquad z^{k_1}\^z^{k_2}\^A(z^{k_3})\^\dots+\dots.
\end{align*}
In the case when $A$ has zero grading, that is, when $A={\rm
Diag}(\dots,a_{-1},a_0,a_1,\dots)$ is diagonal, the application of
the above formula may lead to divergency. Therefore, the action
of~$\h A$ should be regularized and we set, by definition,
$$\widehat A(z^{k_1}\^z^{k_2}\^\dots)=
 \sum_{i=1}^\infty(a_{k_i}-a_{-i})\;z^{k_1}\^z^{k_2}\^\dots.$$

\begin{example}
Consider the operator~$z\pd{}{z}\in\gl$, $z^k\mapsto k\,z^k$. Its
image~$\h{z\,\D\!/\!\D z}$ is called the \emph{energy operator}.
Every basic vector $v_\l$ is an eigenvector for this operator with
the corresponding eigenvalue $|\l|=\sum\l_i$. The eigenvalues of
the energy operator give rise to the grading on the
space~$\L^{\frac\infty2}V$. Under boson-fermion correspondence
this grading corresponds to the quasihomogeneous grading in the
space of power series in the variables $p_i$ with $\deg p_i=i$. In
other words, $\h{z\,\D\!/\!\D z}=\sum_{i=1}^\infty
i\,p_i\,\pd{}{p_i}$.
\end{example}

The correspondence $A\mapsto \h A$ is \emph{not} a Lie algebra
homomorphism. For example, the operators of the multiplication by
$x^m$ commute for different~$m$, but one can easily compute that
$[\h{{x^m}},\h{{x^n}}]=n\,\d_{m,-n}$. In fact, one has $[\h A,\h
B]=\h{[A,B]}+{}$(scalar operator). Therefore, we obtain not a
linear but a \emph{projective} representation of the Lie
algebra~$\gl$. More exactly, the value of the scalar correction
forms a cocycle. Therefore, we obtain a linear representation not
of the algebra~$\gl$ itself but of its one-dimensional central
extension denoted by~$\h\gl$. This algebra is generated by the
operators of the form~$\h A$ and by the scalar operators.

Remark, however, that the most part of the operators used in the
present paper lie in the subalgebra of upper triangular matrices
(that is, the operators of non-negative grading). The cocycle is
trivial on this subalgebra and the correspondence $A\mapsto\h A$
is a Lie algebra homomorphism.

\begin{proposition}
The action of the Lie group of the algebra~$\h\gl$ preserves the
set of decomposable vectors.
\end{proposition}

Indeed, if the matrix of~$A$ has zero diagonal entries, then
$$e^{t \h A}(\f_{1}\^\f_{2}\^\dots)=
 e^{tA}\f_{1}\^e^{tA}\f_{2}\^\dots.$$
 Similarly, if $A={\rm Diag}(\dots,a_{-1},a_0,a_1,\dots)$, then
$$e^{t \h A}(\f_{1}\^\f_{2}\^\dots)=
 e^{t(A-a_1)}\f_{1}\^e^{t(A-a_2)}\f_{2}\^\dots.$$

The proposition suggests a natural way to construct solutions of
the KP hierarchy. It is sufficient to pick any transformation from
the corresponding Lie group and to apply it to the vacuum vector.
The resulting vector corresponds under boson-fermion
correspondence to the $\t$-function of some solution. In order to
apply this procedure in practice, it is useful to have an explicit
description of the action of the Lie algebra~$\h\gl$ in terms of
the variables~$p_i$.

\begin{proposition}\label{bosferm} The correspondence between the differential
operators in~$z$ of order at most two and the action of these
operators in $\C[[p_1,p_2,\dots]]$ is presented in the following
table
$${\renewcommand{\arraystretch}{1.5}
\begin{array}{|c|c|c|}\hline
 \hbox{\small Notation}&\hbox{\small Action in }V=\C[[z]][z^{-1}]&
   \hbox{\small Action in }\C[[p_1,p_2,\dots]]\\\hline
 a_{m},~m>0&z^{m}&p_m\\
 a_{-m},~m>0&z^{-m}&m\,\pd{}{p_m}\\
 a_0&1&0\\
 \L_m&z^m(z\pd{}{z}+\frac{m+1}2)&
    \frac12\sum_{i=-\infty}^\infty {:}a_ia_{m-i}{:}\\
 M_m&z^m(\frac12(z\pd{}{z})^2+\frac{m+1}2z\pd{}{z}+\frac{(m+1)(m+2)}{12})&
    \frac16\sum_{i,j=-\infty}^\infty {:}a_ia_ja_{m-i-j}{:}
 \\\hline
\end{array}}$$
\end{proposition}

The equality $\h{{z^{m}}}=a_{m}=p_m$ can be used for an
independent invariant definition of the boson-fermion
correspondence. Namely, the polynomial (or the formal series)
$P(p_1,p_2,\dots)$ corresponds to the vector
$P(a_{1},a_{2},\dots)v_\varnothing$ of the semi-infinite wedge
space. Conversely, any function can be recovered from its partial
derivatives by the Taylor formula
$$P(p_1,p_2,\dots)=e^{\sum
 p_m\pd{}{q_m}}P(q_1,q_2,\dots)\Bigm|_{q=0}.$$
Therefore the equality $\h{{z^{-m}}}=a_{-m}=m\pd{}{p_m}$ can be
used for the inverse homomorphism of the boson-fermion
correspondence: the vector $v\in\L^{\frac\infty2}V$ corresponds to
the series
$$v\longleftrightarrow \langle e^{\sum \frac{p_m
 a_{-m}}{m}}v\rangle_0,$$
where $\langle\cdot\rangle_0$ denotes the coefficient of the
vacuum vector. If these equalities are taken as the definition of
the boson-fermion correspondence, then its coordinate presentation
given in the beginning of this section can serve as the definition
of Schur functions.

For the operators~$\L_m$ the correspondence of the table is
proved, for example, in~\cite{KS} (in a slightly different
normalization). The case of the cut-and-join operator~$M_0$ is
treated in~\cite{KL} in the relationship with the Hurwitz theory:
this operator is diagonal and the correspondence is established by
comparison of the eigenvalues. Finally, for the operators~$M_n$,
$n\ne0$, the correspondence follows from the commutating relations
$$2\,n\,M_{n}=[M_0,\L_n]-\frac{n^3-n}{12}a_{n},$$
that can be checked independently both in the space of
differential operators in one variable~$z$ and in the space of
differential operators in the variables~$p_1,p_2,\dots$.

\begin{example}
Since $\h1=0$, the cut-and-join operator can be represented in the
form~$M_0=\h{\frac12(z\pd{}{z})^2+\frac12z\pd{}{z}+\frac16}=\h{{\frac12(z\pd{}{z}+\frac12)^2}}$.
Therefore the exponent $e^H$ of the generating function for the
Hurwitz numbers (see Example~\ref{eH1}) corresponds to the
infinite wedge product
$$e^{\b M_0}e^{a_{-1}}v_\varnothing=\f_1\^\f_2\^\f_3\^\dots,$$
where
$$\f_k=
 e^{\frac\b2[(z\pd{}{z}+\frac12)^2-(\frac12-k)^2]}e^{z}z^{-k}
 =\sum_{i=0}^\infty
 e^{\frac\b2[(i-k+\frac12)^2-(\frac12-k)^2]}\frac{z^{-k+i}}{i!}.$$
\end{example}

\section{Once again about the change of coordinates in the
ELSV}\label{transf}

In the framework of the boson-fermion correspondence, the proof of
Theorem~\ref{th5} takes the following form. Consider the
transformation~$\Xi$ sending any function $\Phi(p)$ to the
function $\Psi(q)=(\Phi-H_{0,2})|_{p\to p(q)}$. This
transformation belongs to the Lie group of the Lie
algebra~$\h\gl$, therefore, it can be considered, via the
boson-fermion correspondence, as a linear transformation~$\Xi$ of
the space~$L$ of Laurent series in~$z$.

\begin{proposition}
The transformation~$\Xi$ is given explicitly by
\begin{equation}\label{chPhi}
\Xi:\f(z)\mapsto \psi(\b,z)
 =(1+\b z)^{-\frac32}e^{-\frac{\b z}{2(1+\b z)}}
  \f\Bigl(\tfrac{z}{1+\b z}e^{-\frac{\b z}{1+\b z}}\Bigr).
\end{equation}
\end{proposition}

\emph{Proof}. The transformation $\Xi$ is obtained by integrating
the differential equation of Lemma~\ref{inf}. Using the
correspondence from the table of Proposition~\ref{bosferm} one can
rewrite this equation as the following one:
\begin{equation}\label{eqbeta}
\begin{aligned}
\pd{\psi}{\b}&=-(2\,\L_{1}+\b\,\L_{2})\psi\\
      &=-2\Bigl(z^2\pd{\psi}{z}+z\,\psi\Bigr)
         -\b\Bigl(z^3\pd{\psi}{z}+\frac32z^2\,\psi\Bigr)\\
&=-\Bigl(2\,z^2+\b\,z^3\Bigr)\pd{\psi}{z}
 -\Bigl(2\,z+\frac32\b\,z^2\Bigr)\psi.
\end{aligned}
\end{equation}
This linear PDE can be solved explicitly using, for example, the
method of characteristics. The solution is presented in the
formula of the proposition (as long as the formula is presented,
there is no difficulty to check its validity by the direct
substitution to the equation). The proposition is proved.

\medskip
The statement of the proposition is applicable also if the initial
function~$\f$ depends on additional parameters. Consider the
transformation $e^{\g M_0}e^{a_1}$ sending the vacuum vector to
the exponent $e^{H(\g;p)}$ of the generating function for Hurwitz
numbers (see Example~\ref{eH1}; we changed temporarily from~$\b$
to~$\g$ the notation for the parameter in this transformation).
The Laurent series obtained by the action of this transformation
satisfy the cut-and-join equation
$$
\pd{\f}{\g}=M_0\f
 =\Bigl(\tfrac12\bigl(z\pd{}{z}\bigr)^2
           +\tfrac12z\pd{}{z}+\tfrac1{12}\Bigr)\,\f.$$
The change~\eqref{chPhi} acts on this equation and the direct
substitution shows that for the transformed function~$\psi$ the
equation takes the form
\begin{equation}\label{eqgamma}
\pd{\psi}{\g}=\Bigl(M_0+4\b M_1+6\b^2 M_2+4\b^3 M_3+\b^4 M_4
 +\frac14\b^2z^2+\frac13\b^3z^3+\frac18\b^4z^4\Bigr)\psi
\end{equation}
Since $\g=\b$, we get that the total derivative over the parameter
is given by the sum of the operators on the right hand sides
of~\eqref{eqbeta} and~\eqref{eqgamma}. Taking into account the
change $\g=\b=u^3$ and the corrections arisen from addition of
linear terms $H_{0,1}$ and from the rescaling $q_i\mapsto
u^{-4i}q_i$ that can be written in terms of Laurent series as the
rescaling $z\mapsto u^{-4} z$, we obtain finally the equation of
Theorem~\ref{newcaj}.

All these computations being slightly cumbersome are nevertheless
absolutely elementary; they are fulfilled in the space of
differential operators (of order not greater than~$2$) in one
variable~$z$. In particular, we have finite number of summands at
every step of computation.

\section{Remark on the $\lambda_g$ conjecture}
\def\tT{{\widetilde T}}
\def\Ftop{F^{\rm top}}
\def\tPsi{{\widetilde G}}
\def\Sym{{\rm Sym}}
\let\S\Sigma
The conjecture of C. Faber (proved by now in several
different ways) asserts the equality
\begin{equation}\label{Fab}
\int_{\ocM_{g,n}}\lambda_g\psi_1^{d_1}\dots\psi_n^{d_n}=
 \binom{2g-3+n}{d_1,\dots,d_n}c_g,
\end{equation}
where $\sum d_i=2g-3+n$ and where the constant $c_g$ depends only
on~$g$ but is independent of~$n$ and of the exponents~$d_i$. One
of the simplest proofs of this equality is given in~\cite{GJV2}.
Using the language of the present paper this proof becomes even
more transparent and is reduced to the following.

\begin{theorem}\label{Ftop}
Let $\Ftop$ denote the generating series for the top Hodge
integrals,
$$\Ftop=\sum(-1)^g\langle\l_g\t_0^{k_0}\t_1^{k_1}\dots\rangle_g
 \frac{t_0^{k_0}}{k_0!}\frac{t_1^{k_1}}{k_1!}\dots.$$
Then this series is subject to the equation
 $$\Ftop-\sum_{i\ge0}t_i\pd{\Ftop}{t_i}+\frac12
  \sum_{i,j\ge0}\binom{i+j}{i}t_it_j\pd{\Ftop}{t_{i+i-1}}
  +\frac{t_0^3}{3}=0.$$
\end{theorem}

Denoting by $F^{(m)}$ the homogeneous summand of degree~$m$
in~$\Ftop$ we can rewrite the equation of the theorem in the form
$$F^{(m+1)}=\frac1mA F^{(m)}, \qquad A=
 \frac12\sum_{i,j\ge0}\tbinom{i+j}{i}t_it_j\pd{}{t_{i+i-1}}.$$
This recursive equation allows one to recover the whole series by
induction from its linear part $F^{(1)}=\sum c_gt_{2g-2}$, where
$c_{g}=\langle\l_g\t_{2g-2}\rangle_g$. The fact that this
procedure leads to the equality of Faber's conjecture can be
proved inductively by elementary considerations. Set
$$P_{m,d}=\frac1{m!}\sum_{d_1+\dots+d_m=d}
 \binom{d}{d_1,\dots,d_m}t_{d_1}\dots t_{d_m}.$$
Then the equality~\eqref{Fab} would follow from the relation
$$A P_{m,d}=m\,P_{m+1,d+1}.$$

Consider the linear operation $\S_m$ defined by
$$\S_m \,t_{d_1}\dots t_{d_m}=\Sym_m x_1^{d_1}\dots x_m^{d_m},$$
where $\Sym_mf(x_1,\dots,x_m)=\sum_{\s\in
S(m)}f(x_{\s(1)},\dots,x_{\s(m)})$. The operation $\S_m$ provides
an isomorphism between the space of homogeneous of degree~$m$
polynomials in $t_0,t_1,\dots$ and the space of all symmetric
polynomials in $x_1,\dots,x_m$. Besides, we have, by definition,
$\S_mP_{m,d}=(x_1+\dots+x_m)^d$. Now we compute
$$\S_{m+1}\,A\, t_{d_1}\dots t_{d_m}=
 \frac12\Sym_{m+1}\sum_{k=1}^m
 x_1^{d_1}\dots(x_k+x_{m+1})^{d_k+1}\dots x_m^{d_m}.$$
Therefore,
\begin{align*}
\S_{n+1}AP_{m,d}&
 =\frac1{2\,m!}\Sym_{m+1}\sum_{k=1}^m(x_k+x_{m+1})
  (x_1+\dots+(x_k+x_{m+1})+\dots+x_m)^d\\
 &=\frac1{2\,m!}(x_1+\dots+x_{m+1})^m
    \Sym_{m+1}(x_1+\dots+x_m+m\,x_{m+1})\\
 &=m\,(x_1+\dots+x_{m+1})^{m+1}=m\,\S_{m+1}P_{m+1,d+1}.
\end{align*}
This proves the relation~$A P_{m,d}=m\,P_{m+1,d+1}$, and hence,
the equality~\eqref{Fab} of the $\l_g$-conjecture.

\medskip
For the proof of Theorem~\ref{Ftop} consider the generating
series~$G$ of Theorem~\ref{th2} for Hodge integrals. Set
$v=u^{-3}=\b^{-1}$. A simple counting of dimensions shows that the
series~$G$ can be rewritten in the form
$$ G(q)=v^{-1}\sum_{j,k_0,k_1,\dots}(-1)^{g-j}
 \langle \l_{g-j}\,\t_{0}^{k_0}\,\t_{1}^{k_1}\dots\rangle v^{j}
 \frac{{\tT_0}^{k_0}}{k_0!}\frac{{\tT_1}^{k_1}}{k_1!}\dots,
 \qquad \tT_k=v^{\frac{1-k}3}T_k,$$
where $T_k$ are linear functions in the variables~$q_i$ defined
by~\eqref{D}. Let us make one more change of coordinates passing
to the variables $r_1,r_2,\dots$ related to~$q_k$ by the
equalities
$$q_k=\sum_{i=1}^k\binom{k}{i}(-1)^{k-i}v^{i-\frac k3-1} r_i.$$
Under the identification $q_k\leftrightarrow z^k$ the
variable~$r_k$ corresponds to the polynomial $r_k\leftrightarrow
v\,(v^{-\frac23}z+v^{-1})^k-v^{1-k}$. Up to a rescaling, passing
to the new variables $r_k$ is equivalent to the shift $z\mapsto
z+1$ in the space of polynomials in~$z$. In terms of new variables
the recursive equations~\eqref{D} take the following form
$$\tT_0=r_0,\qquad
\tT_{k+1}=\sum_{m\ge1}m\,(v\,r_{m+2}-r_{m+1})\,\pd{}{r_m}\tT_k.
$$
It follows immediately by induction that the summand of the
smallest degree in~$v$ in~$\tT_k$ is equal to $(-1)^k k! r_{k+1}$:
\begin{align*}
 \tT_0&=~r_1,\\
 \tT_1&=-r_2 + v\,r_3,\\
 \tT_2&=~2\,r_3 - 5\,v\,q_4 + 3\,v^2\,r_5,\\
 \tT_3&=-6\,r_4 + 26\,v\,r_5 - 35\,v^2\,r_6 + 15\,v^3\,r_7,\\
   &\qquad\dots.
\end{align*}

Denote by $\Psi=\Psi(v;r_1,r_2,\dots)$ the series $v\,G$ expressed
in terms of the variables $r_k$. We conclude

\begin{proposition}\label{Psi} The series $\Psi$ is a power series with respect to the
parameter~$v$; for $v=0$ this series turns into the generating
function~$\Ftop$ for the top Hodge integrals after the
substitution
$$t_k=\tT_k|_{v=0}=(-1)^k k! r_{k+1}.$$
\end{proposition}

We would like to apply the change of variables to the
`cat-and-join' equation of Theorem~\ref{newcaj} for~$G$ in order
to obtain the corresponding equation for~$\Psi$. Remark that the
change corresponding to the shift $z\mapsto z+1$ in the space of
polynomials in~$z$ preserves the KP hierarchy (in spite of the
fact that the action of this shift is not defined on the space of
Laurent series). However, the subsequent rescaling is not
quasihomogeneous and destroys the hierarchy. Therefore, it is
convenient to introduce `intermediate' variables $\tilde
r_i=r_i/v$ and to denote by $\tPsi$ the series $G$ expressed in
terms of these variables. Then we have
\begin{equation}\label{renorm}
\Psi(r_1,r_2,\dots)=v\,\tPsi(r_1/v,r_2/v,\dots).
\end{equation}

The transition from the variables $q_i$ to $\tilde r_i$
corresponds to the linear non-homogeneous change $z\mapsto
v^{-\frac23}z+v^{-1}$ in the space of polynomials in~$z$. Applying
this change to the differential operator in~$z$ corresponding to
the right hand side of the equality of Theorem~\ref{newcaj} we
compute the action of this change on the bosonic side as well. As
a result of these computations (that are absolutely elementary
since they do not require consideration of infinite sums) we
arrive at the following statement.

\begin{proposition} The series $\tPsi$ in the variables $\tilde r_k$
is a solution of the KP hierarchy {\rm(}identically in~$v${\rm)}
and its dependence in the parameter~$v$ is described by the
differential equation
$$\pd{e^\tPsi}{v}=\Bigl(-M_2+2\,v\,M_3-M_4+\L_1+\frac{v}{6}\,a_3
 -\frac{v^2}{8}\,a_4\Bigr)e^\tPsi.$$
\end{proposition}

The final rescaling~\eqref{renorm} exits the space of solutions of
the KP hierarchy. Fortunately, this change is simple enough in
order to be able to compute its action on differential operators
in the variables $r_1,r_2,\dots$ directly. We get the following.

\begin{corollary}
The series $\Psi$ of Proposition~{\rm\ref{Psi}} satisfies the
equation
\begin{multline*}
v\pd{\Psi}{v}=\Psi+\sum_{i\ge1}(v\,i\,r_{i+1}-r_i)\pd{\Psi}{r_i}
\\+\frac12\sum_{\substack{i+j=k\\i,j,k\ge1}}\Bigl(
 i\,j\,v\,(-r_{k+2}+2\,v\,r_{k+3}-v^2r_{k+4})
 \bigl(v\pd{^2\Psi}{r_i\D r_j}+\pd{\Psi}{r_i}\pd{\Psi}{r_j}\bigr)\\
 +r_ir_j\bigl(-(k-2)\pd{\Psi}{r_{k-2}}+2\,v\,(k-3)\pd{\Psi}{r_{k-3}}
 -v^2\,(k-4)\pd{\Psi}{r_{k-4}}\bigr)\Bigr)\\
 +\frac13r_1^3-\frac{v}{2}r_1^2r_2+\frac{v^2}{6}r_3-\frac{v^3}{8}r_4.
\end{multline*}
\end{corollary}

Setting~$v=0$ in this equation and denoting $t_k=(-1)^kk!r_{k+1}$
we obtain the required equation of Theorem~\ref{Ftop}.

\section{Reduction of Hodge integrals}\label{secHodge}

For completeness, we review in this section a formula expressing
the Hodge integrals in terms of intersection numbers of just
$\psi$-classes. Denote by $\cF(u,T_0,T_1,\dots)$ the generating
series~\eqref{preHodge} for Hodge integrals written in terms of
$T$-variables. Then $F(T)=\cF(0,T)$ is the Witten's potential for
intersection numbers of $\psi$-classes. The algorithm outlined
in~\cite{F} allowing one to recover $\cF$ from~$F$ is essentially
equivalent to the following formula (cf.~\cite{FP},\cite{G}):
\begin{equation}
\begin{aligned}
 &\qquad e^\cF=e^{W}e^F,\\
  W&=-\sum_{k=1}^\infty
  \tfrac{B_{2k}\,u^{2(2k-1)}}{2k(2k-1)}\Bigr(\pd{}{t_{2k}}-\sum_{i=0}^\infty t_i\pd{}{t_{i+2k-1}}
 +\frac12\sum_{i+j=2k-2}(-1)^i\pd{^2}{t_i\D t_j}\Bigl),
 \end{aligned}
 \label{Gi}\end{equation}
%\begin{equation}
%\begin{aligned}
% e^\cF&=e^{-\sum_{k=1}^\infty
% u^{2(2k-1)}W_{2k-1}}e^F,\\
% W_m&=\tfrac{B_{m+1}}{m(m+1)}\Bigr(\pd{}{t_{m+1}}-\sum_{i=0}^\infty t_i\pd{}{t_{i+m}}
% +\frac12\sum_{i+j=m-1}(-1)^i\pd{^2}{t_i\D t_j}\Bigl),
% \end{aligned}
% \label{Gi}\end{equation}
where $B_m$'s are the Bernoulli numbers.

The Chern classes $\l_i$ of the Hodge bundle over $\ocM_{g,n}$ can
be expressed via the homogeneous components of its Chern character
$\ch=g+\ch_1+\ch_2+\dots$ by the identity
$$1-u^2\l_1+u^4\l_2-u^6\l_3+\dots=e^{-\sum_{m=1}^\infty
(m-1)!\ch_{m} u^{2m}}.$$

The classes $\ch_m$ can be computed from the Mumford's
theorem~\cite{Mu}. It claims that the even components $\ch_{2k}$
vanish, and for odd $m=2k-1$ one has
$$(m{-}1)!\,\ch_m=\frac{B_{m+1}}{m(m+1)}\Bigl(\pi_*(\psi_{n+1}^{m+1})-
 \sum_{i=1}^n\psi_i^m+\frac12 j_*\bigl(\sum\nolimits_{i+j=m-1}
 (-1)^i\psi_{n+1}^i\psi_{n+2}^j\bigr)\Bigr),$$
where $\pi:\ocM_{g,n+1}\to\ocM_{g,n}$ is the natural forgetful map
and $j:\Delta\to\ocM_{g,n}$ is the double cover of the boundary
divisor considered itself as a moduli space of curves with two
additional markings corresponding to the two branches at the
double point, the markings being numbered by $n+1$ and $n+2$,
respectively.

It is straightforward to check that this relation for $\ch_m$
written in terms of intersection numbers is equivalent
to~\eqref{Gi}.

Remark that the operator $W$ participating in~\eqref{Gi} does not
belong to $\h\gl$. Therefore, its usage is not convenient in the
context of the present paper.

\end{document}